\documentclass[12pt]{article}

\usepackage{a4wide}
\usepackage{amssymb}
\usepackage{amsfonts}
\usepackage{amsmath}

\date{}

\newtheorem{proposition}{Proposition}[section]
\newtheorem{theorem}[proposition]{Theorem}
\newtheorem{lemma}[proposition]{Lemma}

\newtheorem{corollary}[proposition]{Corollary}

\def\fd{ {\rm fd}}
\def\lfd{ {\rm lfd}}
\def\fdim{ {\rm d}}

\def\GK{{\rm  GK}\,}
\def\Kdim{{\rm K.dim }\,}

\def\der{\partial }

\def\nFM0{{\nu }_{F,M_0}}
\def\nFN0{{\nu }_{F,N_0}}
\def\nGN0{{\nu }_{G,N_0}}

\def\N0{ {\bf N}_0 }

\def\t{\otimes}
\def\g{\gamma}

\def\ra{\rightarrow}

\def\Xpm{X^{\pm }}

\def\l1{{\lambda}_1}

\def\a{\alpha}
\def\a0{ {\alpha }_0}
\def\a1{ {\alpha }_1}

\def\l{\lambda}
\def\o{\omega}

\def\nFGM0{{\nu }_{F,G,M_0}}

\def\nFN0{{\nu}_{F,N_0}}

\def\sm{{\sigma}^m}

\def\sm1{{\sigma}^{-1}}

\def\smtp1{{\sigma}^{-t+1}}

\def\o{\omega }
\def\S1{S^{-1}}

\def\Xpm1{X^{\pm 1}_1}

\def\sPM1{{\sigma }^{\pm 1}}
\def\sMP1{{\sigma }^{\mp 1 }}

\def\b{\beta}
\def\d{\delta}

\def\L{\Lambda}

\def\OO{{\cal O}}

\def\CD{{\cal D}}

\def\Ytm1{Y^{t-1}}
\def\Yim1{Y^{i-1}}

\def\CM{{\cal M}}

\def\CF{{\cal F}}
\def\CG{{\cal G}}

\def\Der{{\rm Der }}
\def\ad{{\rm ad }}
\def\dim{{\rm dim }}

\def\ker{ {\rm ker } }

\def\gr{ {\rm gr} }

\def\D{ \Delta }
\def\Ev{ {\rm Ev} }

\def\SL2Z{ {\rm SL}_2({\bf Z}) }

\def\Gp1{ G^{1 , 1 } }
\def\P11{ P^{-1 , 1 } }
\def\Pp1{ P^{1 , 1 } }

\def\nCLsr{{}^\nu\kern-2pt {\cal L}^{\sigma , \rho  }}
\def\nP{{}^\nu \kern-2pt P}
\def\nL{{}^\nu\kern-2pt L}
\def\nLL{{}^\nu\kern-2pt \Lambda}
\def\nPsr{{}^\nu\kern-2pt P^{\sigma , \rho  }}
\def\nLsr{{}^\nu\kern-2pt L^{\sigma , \rho  }}
\def\nuCL{{}^\nu\kern-2pt  {\cal L}}
\def\nCLsr{{}^\nu\kern-2pt {\cal L}^{\sigma , \rho  }}
\def\nCL1m{{}^\nu\kern-2pt {\cal L}^{-1 , 1  }}
\def\x1nu{x^\frac{1}{\nu}}
\def\xm1nu{x^{-\frac{1}{\nu}}}

\def\trdeg{{\rm tr.deg}}

\def\CB{{\cal B}}
\def\SL2Z{{\rm SL}_2( \mathbb{Z})}
\def\hSL2Z{\widehat{{\rm SL}}_2( \mathbb{Z})}

\def\bwnm1{ \overline{w}_n^{-1}}

\def\twnm1{ \widetilde{w}_n^{-1}}

\def\o0{\overline{0}}
\def\o1{\overline{1}}

\def\o{\omega}

\def\CB{{\cal B}}
\def\SL2Z{{\rm SL}_2( \mathbb{Z})}
\def\hSL2Z{\widehat{{\rm SL}}_2( \mathbb{Z})}

\def\bwnm1{ \overline{w}_n^{-1}}

\def\twnm1{ \widetilde{w}_n^{-1}}

\def\b0{ \overline{0}}

\def\Gp{\mathfrak{p}}

\def\RA{\Rightarrow}

\def\hol{{\rm hol}}
\def\Cdim{{\rm Cdim}}
\def\CDX{{\cal D}(X)}
\def\GKtrdeg{{\rm GKtr.deg}}
\begin{document}

\author{V.\  Bavula }

\title{Gelfand-Kirillov Dimension of Commutative Subalgebras of Simple Infinite Dimensional
Algebras and their Quotient Division Algebras}

\maketitle

\section{Introduction}

Throughout this paper, $K$ is a  field, a module $M$ over an
algebra $A$  means a {\em left} module denoted ${}_AM$, $\t
=\t_K$.

In contrast to the finite dimensional case, there is no  general
theory of central simple infinite dimensional algebras.  In some
sense, structure of simple {\em finite dimensional} algebras is
`determined' by their {\em maximal commutative subalgebras
(subfields)[see \cite{Pierceb} for example]}. Whether this
statement is true in general is not yet clear. This is certainly
the case for numerous examples of {\em central simple finitely
generated (infinite dimensional) algebras} $A$. A typical example
of $A$ is the ring of differential operators  on a smooth
irreducible affine algebraic variety, its
 coordinate algebra  is a maximal commutative
subalgebra  that completely `determines' the structure of the ring
of differential operators.

{\bf Quantum completely integrable systems}. Let $X$ be a {\em
smooth irreducible affine  algebraic variety} of dimension
$n:=\dim (X)>0$ over a field $K$ of characteristic zero. The {\em
ring of differential operators} $\CD (X)$ is a simple finitely
generated
 $K$-algebra of Gelfand-Kirillov dimension $\GK (\CD (X))=2n$. The
 algebra $\CD (X)$ is a domain and any  commutative finitely
 generated subalgebra in $\CD (X)$ has Krull or Gelfand-Kirillov
 dimension $\leq n$. Recall that
\begin{eqnarray*}
{\rm the \; Gelfand-Kirillov \; dimension}\; \GK (C)& =& {\rm the
\; Krull \; dimension}\; \Kdim (C)\\
 &=&   {\rm the \;  transcendence \; degree} \;
\trdeg_K (C)
\end{eqnarray*}
for every commutative finitely generated algebra $C$ which is a
domain. The algebra  of regular functions $\OO (X)$ on $X$ is a
commutative finitely generated subalgebra of Krull dimension $n$.

{\it Definition}. A {\bf quantum completely integrable system}
(QCIS for short)  is a commutative finitely generated subalgebra
of the algebra of differential operators $\CD (X)$ of Krull
(Gelfand-Kirillov) dimension $n$ (see \cite{BrEtGaQCIS} for
details).

In other words, a QCIS is a commutative finitely generated
subalgebra of $\CD (X)$ of {\em biggest possible} Krull
(Gelfand-Kirillov) dimension. This reformulation  defines a QCIS
for an arbitrary  algebra.

{\it Question}. {\em For a given algebra find an (exact) upper
bound for the Krull (Gelfand-Kirillov) dimension of its
commutative finitely generated subalgebras.}

Surprisingly, it is possible to give such an upper bound only in
terms of `growth', more precisely, in terms of two dimensions (the
Gelfand-Kirillov dimension and the filter dimension)  for any
 central simple finitely generated algebra of finite
Gelfand-Kirillov dimension (Theorem \ref{GKcsuba}) and its
localizations (Theorems \ref{trdsubf}, \ref{c1trds}, and
\ref{c2trds}). Note that the class of central simple finitely
generated algebras of finite Gelfand-Kirillov dimension is a {\em
huge} class of algebras, we are far from understanding structure
of these algebras. Main ingredients of the proofs are the two
filter inequalities (Theorems \ref{FFI} and \ref{SFI}).

For certain classes of algebras and their division algebras the
maximum Gelfand-Kirillov dimension/transcendence degree over the
commutative subalgebras/subfields were found in
\cite{AmitsurPAMS58}, \cite{GK66}, \cite{Mak-Limcom},
\cite{Joseph72HP}, \cite{JosLN74}, \cite{JosgenQ},
\cite{Amitsur-Small78}, and \cite{Resco79}.

{\bf The filter dimension, the first and second filter
inequalities, and Bernstein's inequality}.  Let $A$ be a simple
finitely generated infinite dimensional $K$-algebra. Then
$\dim_K(M)=\infty $ for all nonzero $A$-modules $M$ (the algebra
$A$ is simple, so the $K$-linear map $A\ra {\rm Hom}_K(M,M)$, $a
\mapsto (m \mapsto am)$, is injective, and so $\infty
=\dim_K(A)\leq  \dim_K({\rm Hom}_K(M,M))$ hence
$\dim_K(M)=\infty$). So, the Gelfand-Kirillov dimension (over $K$)
$\GK (M)\geq 1$ for all nonzero $A$-modules $M$.

{\it Definition}. $h_A:=\inf \{ \GK (M)\, | \,  M $  is a nonzero
finitely  generated $A$-module$\}$ is called the {\em holonomic
number} for the algebra $A$.

In \cite{Bavcafd}, the {\bf filter dimension}, $\fd (A)=\fd_K(A)$,
and in \cite{bie98} the {\bf left filter dimension} $\lfd
(A)=\lfd_K(A)$ of simple finitely generated $K$-algebras $A$ were
introduced (see Section \ref{spr1}). In this paper,  $\fdim (A)$
means either the filter dimension $\fd (A)$ or the left filter
dimension $\lfd (A)$ of a simple finitely generated  algebra $A$.
  Both filter dimensions  appear
  naturally when one tries to find a
{\em lower} bound for the holonomic number (Theorem \ref{FFI}) and
an {\em upper} bound (Theorem \ref{SFI}) for the (left and right)
{\em Krull} dimension (in the sense of Rentschler-Gabriel
\cite{Ren-Gab}) of simple finitely generated algebras.

\begin{theorem}\label{FFI}

({\bf The First Filter Inequality},  \cite{Bavcafd, bie98})
 Let $A$ be a simple finitely generated infinite dimensional  algebra.
 Then
$$
 \GK (M)\geq \frac{\GK (A)}{\fdim(A)+\max \{ \fdim(A),
1\} }
$$
 for all nonzero finitely generated $A$-modules $M$ where $\fdim
 =\fd , \lfd$.
\end{theorem}

This theorem is a generalization of {\bf Bernstein's Inequality}
(see Theorem \ref{BerIn}) to a class of simple finitely generated
algebras.

We say that an algebra $A$ is {\em (left) finitely partitive}
 (\cite{MR}, 8.3.17)
if, given any finitely generated $A$-module $M$, there is an
integer $n=n(M)>0$ such that for every strictly descending chain
of $A$-submodules of $M$:
$$M=M_0\supset M_1\supset \cdots \supset M_m$$
with $\GK (M_i/M_{i+1})=\GK (M)$, one has $m\leq n$. McConnell and
Robson write in their book \cite{MR}, 8.3.17,  that ``{\em yet  no
examples are known which fail  to have this property.}''

\begin{theorem}\label{SFI}
({\bf The Second Filter Inequality},  \cite{Bavjafd, bie98})
 Let $A$ be a simple finitely generated finitely partitive algebra with
$\GK (A)<\infty $. Suppose that the Gelfand-Kirillov dimension of
every finitely generated $A$-module is a natural number.
 Then, for any nonzero finitely generated $A$-module $M$,
 the Krull dimension
$$
\Kdim (M)\leq \GK (M) - \frac{\GK (A)}{\fdim (A)+\max \{ \fdim
(A), 1\}}
$$
 where $\fdim  =\fd , \lfd$. In particular,
$$
{\rm K.dim} (A)\leq \GK (A)\left( 1-\frac{1}{\fdim (A)+\max \{
\fdim (A), 1\}} \right).
$$
\end{theorem}
{\it Example}. Let $K$ be a field of characteristic zero, and let
$X$ be a {\em smooth irreducible affine  algebraic variety} of
dimension $n:=\dim (X)>0$. The ring of differential operators $\CD
(X)$ on $X$  is a simple finitely generated infinite dimensional
finitely partitive $K$-algebra with $\GK (\CD (X))=2n$,  $\Kdim
(\CD (X))=n$ \cite{Ren-Gab}, and the Gelfand-Kirillov dimension of
every finitely generated $\CD (X)$-module is a natural number.

\begin{theorem}\label{BerIn}
({\bf Bernstein's Inequality}) $\GK (M)\geq n$ for all nonzero
finitely generated $\CD (X)$-modules $M$.
\end{theorem}

Bernstein \cite{Ber72} proved this inequality for the {\em Weyl
algebra} $A_n=\CD (\mathbb{A}^n)$, the ring of differential
operators on the affine space $\mathbb{A}^n$.

 {\it Definition}. A
nonzero finitely generated $\CD (X)$-module $M$ is called   a {\em
holonomic} module if $\GK (M)=n$ (the least possible
Gelfand-Kirillov dimension).

 This result implies that the holonomic number
$h_{\CD (X)}=n$ since the algebra $\OO (X)$ of regular functions
on $X$ (the coordinate algebra of $X$) is a holonomic $\CD
(X)$-module.

\begin{theorem}\label{fdDX=1}
\cite{Bavjafd, bie98} $\fdim (\CD (X))=1$  where $\fdim
 =\fd , \lfd$.
\end{theorem}
When one puts $\fdim (\CD (X))=1$, $\GK (\CD (X))=2n$, and $\Kdim
(\CD (X))=n$ in the first and second filter inequalities one gets,
in fact, the equalities
$$ n=h_{\CD (X)}\geq \frac{2n}{1+1}=n \;\; {\rm and}\;\; n=\Kdim
(\CD (X))\leq 2n(1-\frac{1}{1+1})=n.$$ There exist other examples
of simple finitely generated infinite dimensional algebras that
are close to the rings of differential operators for which  the
two filter inequalities are also equalities, \cite{Bavcafd} (in
fact, I do not know yet a single example where this is not the
case).

A main goal of this paper is, using the first and the second
filter inequalities, to obtain $(i)$ an {\em upper} bound for the
Gelfand-Kirillov dimension of (maximal) commutative subalgebras of
simple finitely generated infinite dimensional algebras (Theorem
\ref{GKcsuba}), and $(ii)$ an {\em upper} bound for the
transcendence degree of (maximal) subfields of quotient division
rings of (certain) simple finitely generated infinite dimensional
algebras (Theorems \ref{trdsubf} and \ref{c1trds}).

{\bf An upper bound for the Gelfand-Kirillov dimensions of maximal
commutative subalgebras of simple infinite dimensional algebras}.
A $K$-algebra $A$  is called {\em central} if its centre $Z(A)=K$.
\begin{theorem}\label{GKcsuba}
Let $A$ be a central simple finitely generated $K$-algebra of
Gelfand-Kirillov dimension $0<n<\infty $ (over $K$). Let $C$ be a
commutative subalgebra of $A$. Then
$$ \GK (C)\leq \GK (A)\left( 1-\frac{1}{f_A+\max \{
f_A, 1\}} \right)$$ where $f_A:= \max \{ \fdim_{Q_m}(Q_m\t A)\, |
\, 0\leq m\leq  n\}$, $Q_0:=K$, and $Q_m:=K(x_1, \ldots , x_m)$ is
a rational function field in indeterminates $x_1, \ldots , x_m$.
\end{theorem}

A proof of this theorem is given in Section \ref{spr1}. As a
consequence we have a short proof of the following well-known
result.

\begin{corollary}\label{CDXn}
Let $K$ be an algebraically closed  field of characteristic zero,
$X$ be a smooth irreducible affine  algebraic variety of dimension
$n:=\dim (X)>0$, and $C$ be a commutative subalgebra of the ring
of differential operators $\CD (X)$. Then $\GK (C)\leq n$.
\end{corollary}

{\it Proof}. The algebra $\CD (X)$ is central since $K$ is an
algebraically closed field of characteristic zero  \cite{MR}, Ch.
15.  By Theorem \ref{fdDX=1}, $f_{\CD (X)}=1$, and then, by
Theorem \ref{GKcsuba},
$$ \GK (C)\leq 2n (1- \frac{1}{1+1})=n. \; \Box $$

{\it Remark}. For the ring of differential operators $\CD (X)$ the
upper bound of Theorem \ref{GKcsuba} for the Gelfand-Kirillov
dimension  of maximal  commutative subalgebras of $\CD (X)$ is an
{\em exact}  upper  bound since as we mentioned above the  algebra
$\OO (X)$ of regular functions on $X$ is a commutative subalgebra
of $\CD (X)$ of Gelfand-Kirillov dimension $n$.

{\bf An upper bound for the transcendence degree  of maximal
 subfields of quotient division  algebras of simple infinite dimensional algebras}.
 In this paper we prove a general result (Theorem
 \ref{trdsubf}) concerning an upper bound for the transcendence degree  of maximal
 subfields of localizations of (some) simple infinite dimensional
 algebras. Here we only state some of its corollaries which are
 important in applications.

 A $K$-algebra $A$ is said to be a {\em somewhat commutative} if
 it has a finite dimensional  filtration $A=\cup_{i\geq 0}A_i$
 such that the associated graded algebra $\gr (A):=\oplus_{i\geq
 0}A_i/A_{i-1}$ is a commutative finitely generated algebra.
 Typical examples of somewhat commutative algebras are the
 universal enveloping algebra of a finite dimensional Lie algebra
 (and all its factor algebras) and the ring of differential
 operators $\CD (X)$ on  a smooth irreducible affine  algebraic
 variety $X$ over a field of characteristic zero. Every
 somewhat commutative algebra $A$ is a Noetherian  finitely generated finitely
  partitive algebra of finite Gelfand-Kirillov dimension, the Gelfand-Kirillov
  dimension of every finitely generated $A$-modules is an integer, and
 ({\bf Quillen's lemma}): the ring ${\rm End}_A(M)$ is algebraic
 over $K$ (see
  \cite{MR}, Ch. 8 or \cite{KL} for details). If, in addition, the
  algebra $A$ is a domain, then  we denote by $D=D_A$ its {\em
  quotient division  ring (i.e. $D=\S1 A$, $S:=A\backslash \{ 0\}$)}.

\begin{theorem}\label{c1trds}
Let $A$ be a central simple  somewhat commutative infinite
dimensional $K$-algebra which is a domain, and let $D$ be its
quotient division algebra. Let $L$ be a subfield of $D$ that
contains $K$. Then the transcendence degree of the field $L$ (over
$K$)
$$\trdeg_K (L)\leq \GK (A)\left( 1-\frac{1}{f_A+\max \{
f_A, 1\}} \right)$$ where $f_A:= \max \{ \fdim_{Q_m}(Q_m\t A)\, |
\, 0\leq m\leq \GK(A)\}$.
\end{theorem}

\begin{theorem}\label{c2trds}
Let $K$ be an algebraically closed field of characteristic zero,
$\CD (X)$ be the ring of differential operators on a smooth
irreducible affine  algebraic variety $X$ of dimension $n>0$, and
$D(X)$ be the quotient division ring for $\CD (X)$. Let $L$ be a
(commutative) subfield of $D(X)$ that contains $K$. Then
$\trdeg_K(L)\leq n$.
\end{theorem}

{\it Remark}. This inequality is, in fact, an {\em exact} upper
bound for the transcendence degree of subfields in $D(X)$ since
the field of fractions $Q(X)$ for the algebra $\OO (X)$ is a
commutative subfield  of the division ring $D(X)$ with
$\trdeg_K(Q(X))=n$.

Proofs of Theorems \ref{c1trds} and \ref{c2trds} are given in
Section \ref{spr2}.

{\bf An upper bound for the transcendence degree  of maximal
isotropic subalgebras of strongly simple Poisson algebras}. In
Section \ref{maxissPA}, using Theorem \ref{GKcsuba} we prove the
following result

\begin{theorem}\label{IPoGKcsu}
Let $P$ be a strongly simple Poisson algebra, and $C$ be an
isotropic subalgebra of $P$, i.e. $\{ C,C\}=0$. Then
$$ \GK (C)\leq \frac{\GK (A(P))}{2}\left( 1-\frac{1}{f_{A(P)}+\max \{
f_{A(P)}, 1\}} \right)$$ where $f_{A(P)}:= \max \{
\fdim_{Q_m}(Q_m\t A(P))\, | \, 0\leq m\leq  \GK (A(P))\}$.
\end{theorem}

A typical example of the strongly simple Poisson algebra $P$ is
the polynomial algebra $P_{2n}=K[x_1, \ldots , x_{2n}]$ in $2n$
variables over a field $K$ of characteristic zero equipped with
the {\em classical Poisson} bracket (see Section \ref{maxissPA}
for details). Then the algebra $A(P_{2n})$ is the Weyl algebra
$A_{2n}$. Since $\GK (A_{2n})=4n$, $f_{A_{2n}}=1$ we get the
well-known result $$ \GK (C) \leq
\frac{4n}{2}(1-\frac{1}{1+1})=n.$$ This inequality is a sharp one
since the polynomial subalgebra $K[x_1, \ldots, x_n]$ is an
isotropic subalgebra  of $P_{2n}$ of Gelfand-Kirillov dimension
$n$.

{\bf Simple holonomic modules over certain finitely generated
algebras}. In Section \ref{sholmod}, a generalization (Theorem
\ref{simhol}) is given of a construction of A. Braverman, P.
Etingof and D. Gaitsgory (Corollary \ref{csimhol})  that produces
simple holonomic modules (with respect to transcendental field
extensions of the base field).


\section{Proof of Theorem \ref{GKcsuba}}\label{spr1}
{\bf The Gelfand-Kirillov dimension and the filter dimension}. Let
$\CF $ be the set of all functions from the set of natural numbers
$ \mathbb{N}=\{ 0, 1, \ldots \}$ to itself. For each function
$f\in \CF $, the non-negative real number or $\infty $ defined as
$$ \g (f):=\inf \{  r\in \mathbb{R}\, | \, f(i)\leq i^r\; {\rm for
}\; i\gg 0\}$$ is called the  {\em degree} of $f$. The function
$f$ has {\em polynomial growth} if $\g (f)<\infty $. Let $f,g,
p\in \CF $, and $p(i)=p^*(i)$ for $i\gg 0$ where $p^*(t)\in
\mathbb{Q}[t]$ (a polynomial algebra with coefficients from the
field  of rational numbers). Then
\begin{eqnarray*}
\g (f+g)\leq \max \{ \g  (f), \g (g)\}, & & \g (fg)\leq \g (f)+ \g
(g),\\
\g (p)=\deg_t(p^*(t)), & & \g (pg)= \g (p)+ \g
(g).\\
\end{eqnarray*}
Let $A=K\langle a_1, \ldots , a_s\rangle $ be a finitely generated
algebra. The finite dimensional filtration associated with algebra
generators $ a_1, \ldots , a_s$:
$$ A_0:=K\subseteq A_1:=K+\sum_{i=1}^sKa_i\subseteq \cdots \subseteq
A_i:=A_1^i\subseteq \cdots $$ is  called the {\em standard
filtration} for the algebra $A$. Let $M=AM_0$ be a finitely
generated $A$-module where $M_0$ is a finite dimensional
generating subspace. The finite dimensional filtration $\{
M_i:=A_iM_0\}$ is called the {\em standard filtration} for the
$A$-module $M$.

{\it Definition}. $\GK (A):=\g (i\mapsto \dim_K(A_i))$ and $ \GK
(M):=\g (i\mapsto \dim_K(M_i))$ are called the {\bf
Gelfand-Kirillov} dimensions of the algebra $A$ and the $A$-module
$M$ respectively.

It is easy to prove that the Gelfand-Kirillov dimension  of the
algebra (resp. the module)  does not depend on the choice of the
standard filtration of the algebra (resp. and the choice of the
generating subspace of the module).

Suppose, in addition, that the finitely generated algebra $A$ is a
{\em simple} algebra and its centre $Z(A)$ is an {\em algebraic}
field extension of $K$ (the centre of a simple algebra is a
field).  The {\bf return function} $\nu_F \in \CF $ and the {\bf
left return function} $\l_F\in \CF $ for the algebra $A$ with
respect to the standard filtration $F:= \{ A_i\}$ for the algebra
$A$ is defined by the rules:
\begin{eqnarray*}
\nu_F(i)&:=& \min \{ j\in \mathbb{N}\, | \,\,  1\in A_jaA_j\;\,
{\rm
for \; all}\; \,  0\neq a\in A_i\},\\
\l_F(i)&:=& \min \{ j\in \mathbb{N}\, | \, 1\in AaA_j\; \; {\rm
for \; all}\; \; 0\neq a\in A_i\},
\end{eqnarray*}
where $A_jaA_j$ is the vector subspace of the algebra $A$ spanned
over the field $K$ by the elements $xay$ for all $x,y\in A_j$; and
$AaA_j$ is the left ideal of the algebra $A$ generated by the set
$aA_j$. From the definition it is not clear why $\nu_F(i)$ and
 $\l_F(i)$ are finite, the next result proves this.

\begin{lemma}\label{fdfin}
$\l_F(i)\leq \nu_F(i)<\infty $ for $i\geq 0$.
\end{lemma}

{\it Proof}. The first inequality is evident.

 The centre $Z=Z(A)$ of the simple algebra $A$ is a field that
 contains $K$. Let $\{ \o_j\, | \, j\in J\}$ be a $K$-basis for the $K$-vector
space $Z$. Since $\dim_K(A_i)<\infty $, one can find a finitely
many $Z$-linearly independent elements, say $a_1,\ldots , a_s$, of
$A_i$ such that $A_i\subseteq Za_1+\cdots +Za_s$. Next, one can
find a finite subset, say $J'$, of $J$ such that $A_i\subseteq
Va_1+\cdots +Va_s$ where $V=\sum_{j\in J'}K\o_j$. The field $K'$
generated over $K$ by the elements $\o_j$, $j\in J'$, is a finite
field field extension of $K$ (i.e. $\dim_K(K')<\infty $) since
$Z/K$ is algebraic, hence $K'\subseteq A_n$ for some $n\geq 0$.
Clearly, $A_i\subseteq K'a_1+\cdots +K'a_s$.

The $A$-bimodule ${}_AA_A$ is simple with ring of endomorphisms
${\rm End}({}_AA_A)\simeq Z$. By the {\em Density} Theorem,
\cite{Pierceb}, 12.2, for each integer $1\leq j \leq s$, there
exists elements of the algebra $A$, say $x_1^j, \ldots , x_m^j,
y_1^j, \ldots , y_m^j$, $m=m(j)$, such that for all $1\leq l\leq
s$
$$ \sum_{k=1}^m x_k^ja_ly_k^j=\delta_{j,l}, \;\; {\rm the \; Kronecker
\; delta}.$$ Let us fix a natural number, say $d=d_i$, such that
$A_d$ contains all the elements $x_k^j$, $y_k^j$, and the field
$K'$. We claim that $\nu_F(i)\leq 2d$. Let $0\neq a\in A_i$. Then
$a=\l_1a_1+\cdots +\l_sa_s$ for some $\l_i\in K'$. There exists
$\l_j\neq 0$. Then $\sum_{k=1}^m\l_j^{-1}x_k^ja_jy^j_k=1$, and $
 \l_j^{-1}x_k^j, y^j_k\in A_{2d}$. $\Box
$

{\it Definition}. $\fd (A):=\g (i\mapsto \nu_F(i))$ and $\lfd
(A):=\g (i\mapsto \l_F(i))$ are called the {\bf filter dimension}
and the {\bf left filter dimension} of the simple finitely
generated algebra $A$ such that its centre is algebraic over $K$
respectively. By Lemma \ref{fdfin}, $\lfd (A)\leq \fd (A)$.

It is easy to prove that both filter dimensions do not depend on
the choice of the standard filtration $F$, \cite{Bavcafd,bie98}.

{\it Remarks}. 1. If the field $K$ is {\em uncountable} then
automatically the centre $Z(A)$ of a simple finitely generated
algebra $A$ is algebraic over $K$  (since $A$ has a countable
$K$-basis and the rational function field $K(x)$ has uncountable
basis over $K$ since elements $\frac{1}{x+\l }$, $\l \in K$, are
$K$-linearly independent).

2. If a simple finitely generated  algebra $A$ is somewhat
commutative with respect to a filtration $\{ A_i\}$ then the
tensor product of algebras $A\t A^0$ is a somewhat commutative
algebra with respect to the filtration $\{ B_i:=\sum_{j=0}^iA_i\t
A_{i-j}^0\}$ where $A^0$ is the {\em opposite} algebra to $A$. The
algebra $A$ is simple, and so $A$ is a simple $A\t A^0$-module
(i.e. an $A$-bimodule), hence the centre  $Z(A)\simeq {\rm
End}({}_AA_A)$ is algebraic over $K$, by Quillen's lemma.

3. For the definition and properties of the filter dimension of
modules and algebras which are not necessarily simple the reader
is referred to \cite{Bavcafd}.

\begin{proposition}\label{drGK}
Let $A$ and $C$ be  finitely generated algebras such that $C$ is a
commutative domain with field of fractions $Q$, $B:=C\t A$, and
$\CB :=Q\t A$. Let $M$ be a finitely generated $B$-module such
that  $\CM :=\CB \t_BM\neq 0$. Then $ \GK ({}_BM)\geq \GK_Q({}_\CB
\CM )+\GK (C)$.
\end{proposition}

{\it Remark}. $\GK_Q$ stands for the Gelfand-Kirillov dimension
over the field $Q$.

{\it Proof}. Let us fix  standard filtrations  $\{ A_i\}$ and $\{
C_i\}$ for the algebras $A$ and $C$ respectively. Let $h(t)\in
\mathbb{Q}[t]$ be the {\em Hilbert polynomial} for the algebra
$C$, i.e. $\dim_K(C_i)=h(i)$ for $i\gg 0$. Recall that $\GK
(C)=\deg_t(h(t))$. The algebra $B$ has a standard filtration $\{
B_i\}$ which is the tensor product of the standard filtrations
$\{C_i\}$ and $\{ A_i\}$  of the algebras $C$ and $A$, i.e.
$B_i:=\sum_{j=0}^iC_j\t A_{i-j}$. By the assumption, the
$B$-module $M$ is finitely generated, so $M=BM_0$ where $M_0$ is a
finite dimensional generating subspace for $M$. Then the
$B$-module $M$ has a standard filtration $\{ M_i:=B_iM_0\}$. The
$Q$-algebra $\CB $ has a standard (finite dimensional over $Q$)
 filtration  $\{ \CB_i:=Q\t A_i\}$, and the $\CB $-module $\CM $
has a standard (finite dimensional over $Q$) filtration $\{
\CM_i:=\CB_iM_0'=QA_iM_0'\}$ where $M_0'$ is the image of the
vector space $M_0$ under the $B$-module homomorphism $M\ra \CM$,
$m\mapsto m':=1\t_Bm$.

For each $i\geq 0$, one can  fix a $K$-subspace, say $L_i$, of
$A_iM_0'$ such that $\dim_Q(QA_iM_0')=\dim_K(L_i)$. Now,
$B_{2i}\supseteq C_i\t A_i$ implies $\dim_K(B_{2i}M_0)\geq
\dim_K((C_i\t A_i)M_0)$, and $((C_i\t A_i)M_0)'\supseteq C_iL_i$
implies $\dim_K(((C_i\t A_i)M_0)')\geq
\dim_K(C_iL_i)=\dim_K(C_i)\dim_K(L_i)=\dim_K(C_i)\dim_Q(\CM_i)$.
It follows that
\begin{eqnarray*}
\GK ({}_BM)&=&\g (\dim_K(M_i))\geq \g (\dim_K(M_{2i}))=\g
(\dim_K(B_{2i}M_0))\geq \g (\dim_K((C_i\t A_i)M_0))\\
&\geq & \g (\dim_K(((C_i\t A_i)M_0)') \geq \g
(\dim_K(C_i)\dim_Q(\CM_i))\\
&=&\g (\dim_K(C_i))+ \g (\dim_Q(\CM_i))\;\; ({\rm since  }\;\; \g
(\dim_K(C_i))=h(i), \;\; {\rm for }\;\; i\gg 0)\\
&=& \GK (C)+\GK_Q({}_\CB \CM ).\;\; \Box
\end{eqnarray*}

{\bf Proof of Theorem \ref{GKcsuba}}.

Let $P_m=K[x_1, \ldots , x_m]$ be a polynomial algebra over the
field $K$. Then $Q_m$ is its field of fractions and $\GK (P_m)=m$.
Suppose that $P_m$ is a subalgebra of $A$. Then $m=\GK (P_m)\leq
\GK (A)=n$. For each $m\geq 0$, $Q_m\t A$ is a central simple
$Q_m$-algebra (\cite{MR}, 9.6.9) of Gelfand-Kirillov dimension
(over $Q_m$) $\GK_{Q_m}(Q_m\t A)=\GK (A)>0$, hence
$\dim_{Q_m}(Q_m\t A)=\infty$.
\begin{eqnarray*}
\GK (A)&=& \GK({}_AA_A)\geq \GK({}_AA_{P_m})=\GK({}_{ P_m\t A
}A)\;\;\;
(P_m \;\; {\rm is \; commutative})\\
&\geq & \GK_{Q_m}({}_{Q_m\t A}(Q_m\t_{P_m}A))+\GK (P_m) \;\;\;
({\rm
Lemma} \; \ref{drGK})\\
 &\geq & \frac{\GK (A)}{\fdim_{Q_m}(Q_m\t A)+\max \{ \fdim_{Q_m}(Q_m\t A),
 1\}}+m\;\;\; ({\rm Theorem }\; \ref{FFI}).
\end{eqnarray*}
Hence,
$$ m\leq \GK (A)\left( 1-\frac{1}{\fdim_{Q_m} (Q_m\t A)+\max \{
\fdim_{Q_m} (Q_m\t A), 1\}} \right)\leq \GK (A), $$ and so
$$ \GK (C)\leq \GK (A)\left( 1-\frac{1}{f_A+\max \{
f_A, 1\}} \right). \;\; \Box $$


\section{Transcendence Degree of Subfields of the Quotient Division
Algebras of Simple Infinite Dimensional Algebras, Proofs of
Theorems \ref{c1trds} and
\ref{c2trds}}\label{spr2}

Recall that the transcendence degree $\trdeg_K(L)$ of a field
extension $L$ of a field $K$  coincides with the Gelfand-Kirillov
dimension $\GK_K(L)$, and, by {\bf Goldie's Theorem},  a left
Noetherian algebra $A$ which is a domain has a quotient division
ring $D=D_A$ (i.e. $D=\S1 A$ where $S:=A\backslash \{ 0\} $). As a
rule, the division algebra $D$ has {\em infinite} Gelfand-Kirillov
dimension and is {\em not} a finitely generated algebra ({\em eg},
the division ring $D(X)$ of the ring of differential operators
$\CD (X)$ on {\em each} smooth irreducible affine  algebraic
variety $X$ of dimension $n>0$ over a field $K$ of characteristic
zero contains a {\em noncommutative free} subalgebra since
$D(X)\supseteq D(\mathbb{A}^1)$ and the {\em first Weyl division
algebra}  $D(\mathbb{A}^1)$ has this property \cite{Mak-LimFree}).
So, if we want to find an upper bound for the transcendence degree
of subfields in the division ring $D$ we can not apply Theorem
\ref{GKcsuba}. Nevertheless, imposing some natural (mild)
restrictions on the algebra $A$ one can obtain exactly the same
upper bound for the transcendence degree of subfields in the
division ring $D_A$ as the upper bound for the Gelfand-Kirillov
dimension of commutative subalgebras in $A$.

\begin{theorem}\label{trdsubf}
Let $A$ be a simple finitely generated $K$-algebra such that
$0<n:=\GK (A)<\infty $, all the algebras $Q_m\t A$, $m\geq 0$, are
 simple finitely partitive algebras where $Q_0:=K$, $Q_m:=K(x_1,
\ldots , x_m)$ is a rational function field and, for  each $m\geq
0$, the Gelfand-Kirillov dimension (over $Q_m$)  of every finitely
generated $Q_m\t A$-module is a natural number. Let $B=\S1 A$ be
the  localization of the algebra $A$ at a left Ore subset $S$ of
$A$. Let $L$ be a (commutative) subfield of the algebra $B$ that
contains $K$. Then
$$ \trdeg_K(L)\leq \GK (A)\left( 1-\frac{1}{f_A+\max \{
f_A, 1\}} \right)$$ where $f_A:=\max \{ \fdim_{Q_m}(Q_m\t A)\, |
\, 0\leq m \leq  n\}$.
\end{theorem}

{\it Proof}. It follows immediately from  a definition  of the
Gelfand-Kirillov dimension that $\GK_{K'}(K'\t C)=\GK (C)$ for any
$K$-algebra $C$ and any field extension $K'$ of $K$. In
particular, $\GK_{Q_m}(Q_m\t A)=\GK (A)$ for all $m\geq 0$. By
Theorem \ref{SFI},
$$ \Kdim (Q_m\t A)\leq \GK (A)\left( 1-\frac{1}{\fdim_{Q_m} (Q_m\t A)+\max \{
\fdim_{Q_m} (Q_m\t A), 1\}} \right).$$ Let $L$ be a subfield of
the algebra $B$ that contains $K$. Suppose that $L$ contains a
rational function field (isomorphic to) $Q_m$ for some $m\geq 0$.
\begin{eqnarray*}m&=& \trdeg_K(Q_m)\leq \Kdim (Q_m\t Q_m)\\
&\leq & \Kdim (Q_m\t B) \; ({\rm by}\; \cite{MR}, \, 6.5.3\;\;
{\rm since}\;
Q_m\t B \; {\rm is \; a\; free}\;\;  Q_m\t Q_m-{\rm module}) \\
&=& \Kdim (Q_m\t \S1 A)=\Kdim (\S1 (Q_m\t A))\\
&\leq & \Kdim (Q_m\t A)
\; \; ({\rm by}\; \cite{MR}, \, 6.5.3.(ii).(b))\\
&\leq & \GK (A)\left( 1-\frac{1}{\fdim_{Q_m} (Q_m\t A)+\max \{
\fdim_{Q_m} (Q_m\t A), 1\}} \right)\leq \GK (A).
\end{eqnarray*}
Hence
$$\trdeg_K(L)\leq \GK (A)\left( 1-\frac{1}{f_A+\max \{
f_A, 1\}} \right). \;\; \Box $$

{\bf Proof of Theorem \ref{c1trds}}.

The algebra $A$ is a somewhat commutative algebra, so it has a
finite dimensional filtration $A=\cup_{i\geq 0}A_i$ such that the
associated graded algebra is a commutative finitely generated
algebra. For each integer $m\geq 0$, the $Q_m$-algebra $Q_m\t
A=\cup_{i\geq 0}Q_m\t A_i$ has the finite dimensional filtration
 (over $Q_m$) such that the associated graded algebra $\gr (Q_m\t
 A)=\oplus_{i\geq 0}Q_m\t A_i/Q_m\t A_{i-1}\simeq Q_m\t \gr (A)$
 is a commutative finitely generated $Q_m$-algebra. So, $Q_m\t A$
 is a somewhat commutative $Q_m$-algebra.

 By the assumption $\dim_K(A)=\infty$, hence $\dim_K(\gr
 (A))=\infty $ which implies $\GK (\gr (A))>0$, and so $\GK (A)>0$
 (since $\GK (A)=\GK (\gr (A))$). The algebra $A$ is a central
 simple $K$-algebra, so $Q_m\t A$ is a central simple
 $Q_m$-algebra (\cite{MR}, 9.6.9). Now,  Theorem
 \ref{c1trds} follows from Theorem \ref{trdsubf} applied to $B=D$.
 $\Box $

 Let $K$ be a field of characteristic zero, $X$ be a {\em smooth
 irreducible affine  algebraic variety of dimension} $n>0$, $\OO (X)$ be
 its coordinate ring (i.e. the algebra of regular functions on
 $X$).   Recall that the {\em algebra $\CD (X)=\CD (\OO (X))$ of differential
 operators} on $X$  is defined as $\CD (X)=\cup_{i=0}^\infty \,\CD^i
(X)$ where $\CD^0 (X):=\{ u\in {\rm End}_K(\OO (X))\, | \,
ur-ru=0, \; {\rm for \; all}\; r\in \OO (X)\}={\rm End}_{\OO
(X)}(\OO (X))\simeq \OO (X)$, and then inductively
$$ \CD^i (X):=\{ u\in {\rm End}_K(\OO (X))\, | \, ur-ru\in \CD^{i-1} (X),\;
 {\rm for \; all \; }\; r\in \OO (X)\}.$$
Note that the $\{ \CD^i (X)\}$ defines a filtration for the
algebra $\CD (X)$. We say that an element $u\in \CD^i(X)\backslash
\CD^{i-1}(X)$ has order $i$.
\begin{itemize}
\item $\CD (X)$ {\em is a simple somewhat commutative finitely partitive
algebra, a domain.}
\item {\em The algebra $\CD (X)$ is  generated by the algebra $\OO (X)$ and
 the set ${\rm Der}_K (\OO (X))$ of all $K$-derivations of the algebra $\OO
 (X)$}.
\item  {\em The Gelfand-Kirillov dimension} $\GK (\CD (X))=2n$.
\item  {\em The (noncommutative left and right) Krull dimension}
$\Kdim (\CD (X))=n$. \item  $\CD (X)$ {\em is a central algebra
provided $K$ is an algebraically closed field.} \item If $S$ is a
multiplicatively closed subset of $\OO (X)$ then $S$ is an Ore
subset of $\CD (X)$ and $ \CD (\S1 \OO (X))\simeq \S1 \CD (\OO
(X))$  and $ \Der_K(\OO (X))\simeq \S1 \Der_K(\OO (X))$.
\end{itemize}

For proofs of these facts  the reader is referred to [MR], Chapter
15.

{\bf Proof of Theorem \ref{c2trds}}.

Since $Q_m\t \CD_K (\OO (X))\simeq \CD_{Q_m} (Q_m\t \OO (X))$ and
$\fdim (\CD (Q_m\t \OO (X)))=1$ for all $m\geq 0$ we have $f_{\CD
(X)}=1$. Now, Theorem \ref{c2trds} follows from Theorem
\ref{c1trds}, $$ \trdeg_K(L)\leq 2n(1-\frac{1}{1+1})=n. \;  \;
\Box
$$

Following \cite{JosgenQ} for a $K$-algebra $A$ define the {\bf
commutative dimension}
$$ \Cdim (A):=\sup \{ \GK (C)\, | \; C \;\; {\rm is \; a \;
commutative\; subalgebra \; of }\; A\}.$$ The commutative
dimension $\Cdim (A)$ is the largest non-negative integer $m$ such
that the algebra $A$ contains a polynomial algebra in $m$
variables (\cite{JosgenQ}, 1.1, or \cite{MR}, 8.2.14). So, $\Cdim
(A) =\mathbb{N}\cup \{ \infty \}$. If $A$ is a subalgebra of $B$
then $\Cdim (A)\leq \Cdim (B)$.

\begin{corollary}\label{DXYnmnot}
Let $X$ and $Y$ be smooth irreducible affine  algebraic varieties
of dimensions $n$ and $m$ respectively, let $D(X)$ and $D(Y)$ be
quotient division rings for the rings of differential operators
$\CD (X)$ and $\CD (Y)$. Then there is no $K$-algebra embedding
$D(X)\ra D(Y)$ for $n>m$.
\end{corollary}

{\it Proof}. By Theorem \ref{c2trds}, $\Cdim (D(X))=n$ and $\Cdim
(D(Y))=m$. Suppose that there is a $K$-algebra embedding $D(X)\ra
D(Y)$. Then $n=\Cdim (D(X))\leq \Cdim (D(Y))=m$.
 $\Box $

For the Weyl algebras $A_n=\CD ( \mathbb{A}^n)$ and  $A_m=\CD (
\mathbb{A}^m)$ the result above was proved by Gelfand and Kirillov
in \cite{GK66}. They introduced a new  invariant of an  algebra
$A$, so-called the {\em (Gelfand-Kirillov) transcendence degree}
$\GKtrdeg (A)$, and proved that $\GKtrdeg (D_n)=2n$. Recall that
$$ \GKtrdeg (A):=\sup_{V} \inf_{b} \, \GK (K[bV])$$
where $V$ ranges over the finite dimensional subspaces of $A$ and
$b$ ranges over the regular elements of $A$. Another proofs based
on different  ideas were given by A. Joseph \cite{JosLN74} and R.
Resco \cite{Resco79}, see also \cite{MR}, 6.6.19. Joseph's proof
is based on the fact that the centralizer of any isomorphic copy
of the Weyl algebra $A_n$ in its division algebra $D_n:=D(
\mathbb{A}^n)$ reduces to scalars (\cite{JosgenQ}, 4.2), Resco
proved that $\Cdim (D_n)=n$ (\cite{Resco79}, 4.2) using the result
of Rentschler and Gabriel
 \cite{Ren-Gab} that $\Kdim (A_n)=n$  (over an arbitrary field of characteristic
zero).

The next result is a  generalization of Quillen's lemma and is due
to Joseph and Rentschler in \cite{JosgenQ}.

\begin{theorem}\label{JRenQ}
Let $M$ be a finitely generated module  over  a somewhat
commutative algebra $A$. Then $\Cdim ({\rm End}_A(M))\leq \Kdim
(M)$.
\end{theorem}

The next result is due to L. Makar-Limanov.

\begin{theorem}\label{MLcomth}
\cite{Mak-Limcom}. Let $X$ be a smooth irreducible affine
algebraic variety of dimension $n>0$, and let $C$ be a commutative
subalgebra of $\CD (X)$ of Gelfand-Kirillov dimension $n$. Then
its centralizer $C(C, \CD (X))$ is a commutative algebra.
\end{theorem}

As a direct consequence of the previous result we obtain a
characterization of maximal  commutative subalgebras of
Gelfand-Kirillov dimension $n$ in $\CD (X)$.

\begin{lemma}\label{Cmaxn}
Let $X$ be a smooth irreducible affine  algebraic variety of
dimension $n>0$, and let $C$ be a commutative subalgebra of $\CD
(X)$.
 The following statements are equivalent.
\begin{enumerate}
\item $C$ is a maximal  commutative subalgebra of $\CD (X)$ with
$\GK (C)=n$.
\item $C$ is the centralizer in $\CD (X)$ of $n$ commuting
algebraically independent elements of $\CD (X)$.
\item $\GK (C)=n$ and $C$ is the centralizer in $\CD (X)$ of every $n$ commuting
algebraically independent elements of $C$.
\end{enumerate}
\end{lemma}

{\it Proof}. $(1\RA 3)$ Let $T$ be  a subset of $C$ that consists
of $n$ (commuting) algebraically  independent elements. By Theorem
\ref{MLcomth}, the centralizer $C(T)$ of the  set $T$ in $\CD (X)$
is a commutative algebra  that  contains $C$. Therefore, $C(T)=C$
since $C$ is a maximal commutative subalgebra.

$(3\RA 2)$ This implication is evident.

$(2\RA 1)$ Let $C$ be as in the second statement, and $C'$ be a
commutative algebra that contains $C$. Then $C'\subseteq C$ since
$C$ is a centralizer. Therefore, $C$ is a maximal commutative
subalgebra  with Gelfand-Kirillov dimension $n$.  $\Box $

\begin{corollary}\label{coCmaxn}
Let $X$ be a smooth irreducible affine  algebraic variety of
dimension $n>0$, let $C$ and $C'$ be maximal commutative
subalgebras of $\CD (X)$ of Gelfand-Kirillov dimension $n$. Then
either $C=C'$ or, otherwise, $\GK (C\cap C')<n$.
\end{corollary}

{\it Proof}. Suppose that $\GK (C\cap C')=n$. Then one can  choose
a subset, say $T$, of $C\cap C'$ that consists of $n$ (commuting)
algebraically independent elements. By Lemma \ref{Cmaxn}.(3),
$C=C(T, \CD (X))=C'$.  $\Box $

{\it Example}. The polynomial algebras $C=K[x_1, \ldots , x_n]$
and $C'=K[x_1, \ldots , x_m, \der_{m+1}, \ldots , \der_n ]$ are
maximal commutative subalgebras of the Weyl algebra $A_n$ with
$C\cap C'=K[x_1, \ldots , x_m]$. So, the number $m=\GK (C\cap C')$
in Corollary \ref{coCmaxn} can be any natural number between $0$
and $n$.

Let $M$ be a module over a polynomial algebra $K[t]$ where $K$ is
an algebraically closed field (for simplicity). The element $t$ is
called a {\em locally finite element} if $\dim_K(K[t]m)<\infty $
for all $m\in M$, $t$ is a {\em locally nilpotent element} if, for
each $m\in M$, $t^im=0$ for all $i\gg 0$, $t$ is a {\em locally
semi simple element} if $M$ is a semi-simple $K[t]$-module.

Let $T=\{ t_1, \ldots , t_n\}\subseteq \CD (X)$ be a set of {\em
commuting algebraically independent elements}. Let $C(T)=C(T, \CD
(X))=\{ a\in \CD (X)\, | \, at_i=t_ia, \; i=1, \ldots ,
n\}=\cap_{i=1}^n\ker (\ad (t_i))$ be the {\em centralizer} of the
set $T$ in $\CD (X)$. By Theorem \ref{MLcomth}, $C(T)$ is a
commutative algebra. For the set $T$, let $F(T)$ (resp. $N(T)$,
$D(T)$) be the largest subalgebra of $\CD (X)$  on which each
inner derivation $\ad (t_i)$, $i=1,\ldots , n$, is {\em locally
finite} (resp. {\em locally nilpotent}, {\em locally
semi-simple}). Clearly, $C(T)=N(T)\cap D(T)$, $N(T)\subseteq
F(T)$, and $D(T)\subseteq F(T)$. If the field $K$  is
algebraically closed then $$D(T)=\bigoplus_{\l \in \Ev (T)}D(T, \l
)\;\; {\rm and }\;\; F(T)=\bigoplus_{\l \in \Ev (T)}F(T, \l ),
$$ where $\Ev (T):=\{ \l =(\l_1, \ldots ,
\l_n)\in K^n\, | \, [t_i, a]=\l_ia$ for some $0\neq a\in \CD (X)$,
$i=1, \ldots , n\}$ is the set of eigenvalues or weights for $T$,
$D(T, \l ):=\{ a\in \CD(X)\, | \, [t_i, a]=\l_ia, i=1, \ldots ,
n\}$, $F(T, \l ):=\{ a\in \CD(X)\, | \, (\ad (t_i)-\l_i)^{m_i} (
a)=0$ for some $m_1, \ldots , m_n\in \mathbb{N}\}$. $D(T,0)=C(T)$
and $D(T, \l )D(T, \mu )\subseteq D(T, \l +\mu )$ for all $\l ,
\mu \in \Ev (T)$. So, $\Ev (T)$ is an additive sub-semigroup  of
$K^n$ since $\CD (X)$ is a domain.

Similarly, for any set $\D =\{ \d_1, \ldots, \d_t\}$ of commuting
$K$-derivations of the algebra $A$ one can defined the algebras
$C(\D , A)$, $N(\D , A)$, $D(\D , A)$, and $F(\D , A)$.

\begin{lemma}\label{FSS14}
Let $X$ be a smooth irreducible affine algebraic  variety of
dimension $n>0$, $T=\{t_1, \ldots , t_n\}\subseteq \CD (X)$ be a
set of commuting algebraically independent elements. The sets
$S:=K[T]\backslash  \{ 0\}$ and $S_1:=C(T, \CD (X))\backslash \{
0\}$ are Ore subsets of the algebras $C(T, \CD (X))$,  $N(T, \CD
(X))$, $D(T, \CD (X))$, and $F(T, \CD (X))$.
\begin{enumerate}
\item $C(T, D (X))=\S1 C(T, \CD (X))=\S1_1C(T, \CD (X))$.
\item $N(T, D (X))=\S1 N(T, \CD (X))=\S1_1N(T, \CD (X))$.
\item $D(T, D (X))=\S1 D(T, \CD (X))=\S1_1D(T, \CD (X))$, and $\Ev
(T, \CD (X))=\Ev (T, D(X))$ is an additive subgroup of $
\mathbb{Q}^k$, $k\leq n$.
\item $F(T, D (X))=\S1 F(T, \CD (X))=\S1_1F(T, \CD (X))$.
\end{enumerate}
\end{lemma}

{\it Proof}. $4$. It suffices to prove that {\em an arbitrary
element $a$ of the algebra $F':=F(T, D(X))$ has the form $s^{-1}b$
for some $s\in S$ and $b\in \CD (X)$}, since then $b\in F:=F(T,
\CD (X))$, $S$ and $S_1$ are left Ore sets of $F$ (by symmetry,
$S$ and $S_1$ are also right Ore subsets of $F$).

The division algebra $D(X)$ is a module over the polynomial
algebra $K[T]$ where  the action is given by the rule: $t_i*u:=\ad
(t_i)(u)$. The vector space $V=K[T]*a$ has finite  dimension over
$K$ since   $a\in F'$. Therefore, $I:=\{  c\in \CD (X)\, | \,
cV\subseteq \CD (X)\}$ is a nonzero left ideal in $\CD (X)$. The
normalizer $N(I)=\{ c\in \CD (X)\, | \, Ic\subseteq  I\}$ of $I$
in $\CD (X)$ contains $K[T]$ as follows from $It_iV\subseteq
I[t_i, V]+IVt_i\subseteq \CD (X)$. The opposite  algebra
$(N(I)/I)^0$ to  the factor algebra $N(I)/I$ can be canonically
identified with the endomorphism  algebra ${\rm End}_{\CD (X)}(\CD
(X)/I)$ $((N(I)/I)^0\ra {\rm End}_{\CD (X)}(\CD (X)/I)$, $u\mapsto
(c+I\mapsto cu+I))$. Recall that  the opposite algebra $A^0$ to an
algebra $A$ has the same  additive structure as $A$ and
multiplication is defined as $x\cdot y =yx$. Since $\CD (X)$ is a
domain and $I\neq 0$, $\Kdim ( \CD (X)/I)<\Kdim (\CD (X))=n$. By
Theorem \ref{JRenQ},
$$\Cdim ((N(I)/I)^0)\leq \Kdim (\CD (X)/I)<n,$$
hence $K[T]\cap I\neq 0$ since $\GK (K[T])=n$. Take any $0\neq
s\in K[T]\cap I$, then $b:=sa\in \CD (X)$, as required.

$1$ and $2$. Given $s\in S$ and $b\in \CDX$. Then $s^{-1}b\in C(T,
D(X))$ (resp. $s^{-1}b \in N(T, D(X)))$ iff $b\in C(T, \CDX )$
(resp. $b\in N(T, \CDX ))$ and the result follows.

$3$. Statement 4 implies $D(T, D(X))=\S1 D(T, \CDX )=\S1_1D(T,
\CDX )$. Given $\l \in \Ev (T, \CDX )$  and $0\neq a\in D(T, \l ,
\CDX )$. Then $a^{-1}\in D(T, -\l   , \CDX )$ and $sa^{-1}\in \CDX
$ for some $s\in S$. Clearly, $sa^{-1}\in D(T, -\l  , \CDX )$.
Hence $\Ev (T, \CDX )$ is an  additive subgroup in $K^n$ that
coincides with $\Ev (T, D(X))$  since $D(T, D(X))=\S1 D(T, \CDX
)$. Let $\l^1, \ldots , \l^m$, be $ \mathbb{Q}$-linearly
independent elements of $\Ev (T, \CDX )$.  For each $i=1, \ldots ,
m$,  choose $0\neq a_i\in D(T, \l^i, \CDX )$. Using    the $\Ev
(T)$-graded structure of the algebra $D(T, \CDX )$, we see that
the algebra generated by  $T, a_1, \ldots , a_m$ is a polynomial
algebra in $n+m$ variables, so $n+m\leq \GK (\CDX )=2n$ implies
$m\leq n$. $\Box $

The proof of Lemma \ref{FSS14} is based on two facts: a
generalization of Quillen's Lemma (Theorem \ref{JRenQ}) and $\Kdim
(\CD (X))=\Cdim (\CD (X))$. So, repeating word for word this proof
we have a slightly more general result.

\begin{lemma}\label{cFFSS14}
Let a domain $A$ be a somewhat commutative algebra with $n:=\Kdim
(A)=\Cdim (A)$, let $D=D_A$ be its quotient division algebra, and
$T=\{ t_1, \ldots , t_n\}\subseteq A$ be a subset of commuting
algebraically independent elements. Then the results of Lemma
\ref{FSS14} hold with $k\leq \GK (A)-n$, $S$ and $S_1$ are left
Ore subsets of the algebras from  Lemma \ref{FSS14}.
\end{lemma}

{\it Example}. Let $A=U(\CG )$ be the universal enveloping algebra
of a finite dimensional Lie algebra $\CG $ over the field $
\mathbb{C}$ of complex numbers such that $\Kdim (A)=\Cdim (A)$
({\it eg}, $Usl(2)$ since $\Kdim (Usl(2))=2=\Cdim (Usl(2))$).

\begin{corollary}\label{QCmaxDX}
Let $X$ be a smooth irreducible affine algebraic  variety of
dimension $n>0$, and $C$ be a maximal commutative subalgebra  of
$\CDX $ with $\GK (C)=n$. Then its field of fractions $Q(C)$ is a
maximal commutative subfield of the division algebra $D(X)$.
\end{corollary}

{\it Proof}. By Lemma \ref{Cmaxn}, $C=C(T, \CDX )$ for a subset
$T$ of $C$ that consists of  $n$ algebraically independent
elements. Given a subfield $L$ of the division algebra $D(X)$
containing  $Q(C)$. Then $L\subseteq  C(T, D(X))=Q(C)$, by Lemma
\ref{FSS14}.(1). So, $Q(C)$ is a maximal subfield in $D(X)$.
$\Box $

\begin{corollary}\label{OXmax}
Let $X$ be a smooth irreducible affine algebraic  variety of
dimension $n>0$.
\begin{enumerate}
\item The algebra $\OO (X)$ of regular functions on $X$ is a
maximal  commutative subalgebra in $\CDX $ that coincides with its
centralizer $C(\OO (X), \CDX )$.
\item The field of fractions $Q(X)$ of the algebra $\OO (X)$ is a
maximal commutative subfield in the division algebra $D(X)$.
\end{enumerate}
\end{corollary}

{\it Proof}. $1$. By \cite{MR}, 15.2.6, there exists a nonzero
element $s\in \OO (X)$ such that
$$\CD (\OO (X)_s)=\OO
(X)_s[\der_1, \ldots , \der_n]\supseteq A_n:=K\langle x_1, \ldots
, x_n, \der_1, \ldots , \der_n\rangle , \; {\rm the \; Weyl  \;
algebra}, $$ where $\OO (X)_s$ is a localization of the algebra
$\OO (X)$ at the powers of the element $s$; $x_1, \ldots , x_n$
are algebraically independent elements of $\OO (X)_s$; $\der_1,\
\ldots , \der_n$ are commuting $K$-derivations of the algebra $\OO
(X)_s$ satisfying $\der_i(x_j)=\d_{i,j}$, the Kronecker delta. So,
the algebra $\CD (\OO (X)_s)$ contains the Weyl algebra $A_n$, and
the inclusion $A_n=\CD (\mathbb{A}^n)\subseteq \CD (\OO (X)_s)$
respects the canonical filtrations (by the total degree of
derivations).

Let $0\neq c\in C(\OO (X), \CD (X))$ be an element of order $i$.
We have to prove that $i=0$. Suppose to the  contrary that $i>0$.
Then $c=\sum_{ \{ \alpha \in \mathbb{N}^n: \, | \alpha | =i\} }
\l_\alpha  \der^\alpha +\cdots $ where the three dots denote terms
of smaller order, $\alpha =(\alpha_1, \ldots , \alpha_n)$, $|
\alpha | =\alpha_1+\cdots +\alpha_n$, $\der^\alpha=
\der_1^{\alpha_1}\cdots \der_n^{\alpha_n}$. There exists $\alpha $
such that $\l_\alpha \neq 0$. Then $0=\prod_{i=1}^n \ad
(x_i)^{\alpha_i}(c)=(-1)^{| \alpha |}\alpha_1!\cdots
\alpha_n!\l_\alpha \neq  0 $, a contradiction. Therefore, $i=0$.
This implies that $\OO (X)$ is a maximal commutative subalgebra,
by Lemma \ref{Cmaxn}.

$2$. By the first statement and Corollary \ref{QCmaxDX}, $Q(X)$ is
a maximal subfield in $D(X)$. $\Box $

\begin{lemma}\label{GCgr3}
Let $G$ be a semigroup with identity $e$ such that $xy=e$ implies
$yx=e$ for $x,y\in G$. Let a $K$-algebra $B$ be a domain with
$n:=\GK (B)<\infty $. Suppose that the algebra $B$ contains a
simple subalgebra $A$ with $\GK (A)=n$. Then
\begin{enumerate}
\item $B$ is a simple algebra. \item Suppose that
$B=\bigoplus_{g\in G}B_g$ is a $G$-graded algebra and $B_g\neq 0$
for all $g\in G$. Then $G$ is a group. \item Suppose that
$C=\bigoplus_{g\in G}C_g$ is a simple $G$-graded algebra of finite
Gelfand-Kirillov dimension  which is a domain and $C_g\neq 0$ for
all $g\in G$. Then $G$ is a group.
\end{enumerate}
\end{lemma}

{\it Proof}. $1$. Let $I$ be a nonzero ideal of the algebra $B$.
By \cite{MR}, 8.3.5, $\GK (B/I)<\GK (B)$ since $B$ is a domain,
hence $A\cap I\neq 0$ (since otherwise the natural map $A\ra B/I$
were an algebra monomorphism and we would have $n=\GK (A)\leq \GK
(B/I)<n$, a contradiction). The algebra $A$ is simple, so $I\cap
A=A$, hence $I=B$. This proves that $B$ is a simple algebra.

$2$. Since $xy=e$ implies $yx=e$ in $G$ the semigroup $G$ is a
group iff $GgG=G$ for all $g\in G$. Suppose that $G$ is not a
group then $GgG\neq G$ for some element $g\in G$. Then the set
$BB_gB\subseteq \bigoplus_{h\in GgG}B_h$ is a proper ideal in $B$
which contradicts to simplicity of the algebra $B$.

$3$. This is a particular case of statement 2 when $A=B=C$. $\Box
$

\begin{corollary}\label{cGCgr3}
Let a domain $A$ be a simple finitely generated algebra oven an
algebraically closed field $K$ of characteristic zero with $\GK
(A)<\infty $, and let $\D =\{ \d_1, \ldots , \d_t\}$ be a set of
locally finite commuting $K$-derivations of the algebra $A$. Then
$\Ev (\D )\simeq \mathbb{Z}^k$ is a free finitely generated
abelian group of rank $k$ and $k\leq \GK (A)-\GK (C(\D ))$.
\end{corollary}

{\it Proof}. The set $E:=\Ev (\D )$ is an additive sub-semigroup
of $K^t$  since $A$ is a domain. The algebra $A=\bigoplus_{\l \in
E}F(\D , \l )$ is an $E$-graded algebra with $F(\D , \l )\neq 0$
for all $\l \in E$. By Lemma \ref{GCgr3}.(3), $E$ is a subgroup of
$K^t$ since $A$ is a simple domain. The algebra $A$ is finitely
generated, so $E$ is a finitely generated torsion free $
\mathbb{Z}$-module. Hence $E\simeq \mathbb{Z}^k$ for some $k\geq
0$.

Let $\l^1, \ldots , \l^k$ be free generators for the $
\mathbb{Z}$-module $E$, and let $0\neq x_i\in D(\D , \l^i)$ for
each $i$. The algebra $A=\bigoplus_{\l \in E}D(\D , \l )$ is an
$E$-graded domain. So, the left $C(\D )$-submodule
$\bigoplus_{m\in \mathbb{N}^k}C(\D )x^m$ of $A$ is free with the
set  $\{ x^m=x_1^{m_1}\cdots x_k^{m_k}\, | \, m\in \mathbb{N}^k\}$
of free generators. This implies that $\GK (C(\D ))+k\leq \GK
(B)\leq \GK (A)$ where $B$ is the subalgebra of $A$ generated by
$C(\D )$ and $x^m$, $m\in \mathbb{N}^k$.  $\Box $

Let $\d $ be a locally  finite $K$-derivation of an algebra $A$
over an algebraically  closed field $K$ of characteristic zero
(for simplicity). Then $\d $ is a {\em unique} sum $\d =\d_n+\d_s$
of {\em commuting locally nilpotent derivation} $\d_n$ and a {\em
locally semi-simple derivation} $\d_s$. The derivation $\d_s$ is
defined as follows: $\d_s(u)=\l u$ for all $u\in F(\d , \l )$ and
$\l \in \Ev (\d )$. Then $\d_n:=\d -\d_s$. This decomposition is
called the {\bf Jordan decomposition}  for the locally finite
derivation $\d $. Given another locally finite derivation $\d'$ of
the algebra $A$ with Jordan decomposition $\d'=\d_n'+\d_s'$. It is
obvious that {\em the derivations $\d $ and $\d'$ commute iff all
the derivations $\d_n$, $\d_s$, $\d_n'$ , and $\d_s'$ commute}.
{\it Proof}. $(\Rightarrow )$ Suppose that $\d \d' =\d' \d$. Take
$a\in A$, then $V:=K[\d ,\d']a$ is a finite dimensional subspace
of $A$, hence is invariant under the natural action of the
derivations $\d_s $ and $\d_s'$. Clearly, the restrictions of the
derivations $\d_s$ and $\d_s'$ to $V$ are the semi-simple parts of
the restrictions of $\d$ and $\d'$ to $V$ respectively. Since the
restrictions  $\d_s|_V$ and $\d_s'|_V$ are polynomials of  $\d|_V$
and $\d'|_V$ respectively, they commute. So, $\d_s$ and $\d_s'$
commute and then $\d_n$ and $\d_n'$ commute. $\Box$.

{\it Example}. $\d=\sum_{i=1}^m\l_i \frac{\der}{\der
x_i}+\sum_{j=m+1}^n\l_jx_j \frac{\der}{\der x_j}$ is a locally
finite derivation of the polynomial algebra $K[x_1, \ldots ,
x_n]$, and $\d =\d_n +\d_s$, $\d_n =\sum_{i=1}^m\l_i
\frac{\der}{\der x_i}$, $\d_s=\sum_{j=m+1}^n\l_jx_j
\frac{\der}{\der x_j}$, is its Jordan decomposition where
$\l_1,\ldots , \l_n\in K$.

We say that an element $a\in A$ is {\em locally finite} (resp.
{\em locally nilpotent}, {\em locally semi-simple}) if so is the
inner derivation $\ad (a)$. Suppose that {\em all $K$-derivations
of the algebra $A$ are inner. Then every locally finite element
$a$ of $A$ is a sum $a=a_n+a_s$ of a locally nilpotent element
$a_n$ and a locally semi-simple element $a_s$ and they commute. If
$a=a_n'+a_s'$ is another such a sum then $a_n'=a_n+z$ and
$a_s'=a_s-z$ for a unique central element $z\in Z(A)$, and vice
versa}.
 {\it Proof}. Let $\ad (a)=\d_n+\d_s$ be a Jordan decomposition
 for $\ad (a)$. All $K$-derivations of the algebra $A$ are inner,
 so $\d_n=\ad (a_n)$ and $\d_s=\ad (a_s)$ where $a_n\in A$ is a
 locally nilpotent element and $a_s\in A$ is a locally semi-simple
 element. $0=[\ad (a_n), \ad (a_s)]=\ad ([a_n, a_s])$ implies $\l
 :=[a_n, a_s]\in Z(A)$. Since the element $a_s$ is locally
 semi-simple, $\l =0$. Inner derivations $\ad (x)$ and  $\ad (y)$ of
 the algebra $A$ are equal iff
 $x=y+z$ for some $z\in Z(A)$. $\ad (a_n)=\d_n=\ad (a_n')$,
 $\ad (a_s)=\d_s=\ad (a_s')$, $a=a_n+a_s=a_n'+a_s'$, imply
 $a_n'=a_n+z$ and $a_s'=a_s-z$ for a unique $z\in Z(A)$, and vice
 versa. $\Box $

 In particular, we have proved that {\em given a locally
 semi-simple element  $a$, then elements $a$ and $b$ commute iff the
 inner derivations $\ad (a)$ and $\ad (b)$ commute}.

{\it Definition}. For the locally finite element $a$,  the
decomposition $a=a_n+a_s$ above will  be called  a {\bf Jordan
decomposition} for $a$ (it is unique up to an element of the
centre $Z(A)$ as above).

{\it Example}. All the $K$-derivations of the Weyl algebra
$A_n=K\langle x_1, \ldots , x_n ,\der_1, \ldots , \der_n\rangle$
are inner \cite{DixUEA}, 4.6.8. $a=a_n+a_s$,  $a_n =
\sum_{i=1}^m\l_ix_i$, $a_s=\sum_{j=m+1}^n\l_jx_j\der_j$, is the
Jordan decomposition for a locally finite element $a$ where $\l_1,
\ldots , \l_n \in K$.

Let $a$ and $b$ be locally finite elements of the algebra $A$, and
 let $a=a_n+a_s$ and $b=b_n+b_s$ be their Jordan decompositions.
 Then {\em the elements $a$ and $b$ commute iff all the elements
 $a_n$, $a_s$, $b_n$, and $b_s$ commute}. {\it Proof}. Suppose
 that the elements $a$ and $b$ commute then the inner derivations
 $\ad (a) $ and $\ad (b)$ commute, then all the derivations
 $\ad (a_n) $, $\ad (a_s)$, $\ad (b_n) $ and $\ad (b_s)$ commute.
 The elements $a_s$ and $b_s$ are locally semi-simple, hence
 $a_s$ (resp. $b_s$) commute with $b_n$ and $b_s$ (resp. $a_n$ and
 $a_s$). So, all the elements $a_n$, $a_s$, $b_n$, and $b_s$
 commute. The inverse implication is obvious. $\Box $

\begin{corollary}\label{ckGCgr3}
Let  a domain  $A$ be a simple finitely generated algebra over an
algebraically closed field $K$ of characteristic zero such that
every $K$-derivation of the algebra $A$ is inner and $n:=\GK
(A)<\infty $. Let $\D = \{ \d_1, \ldots , \d_t\}$ be a set of
commuting locally finite $K$-derivations of the algebra $A$. Then
\begin{enumerate}
\item $\Ev (\D )\simeq \mathbb{Z}^k$ with $\GK (K\langle \d_{1,
s}, \ldots , \d_{t, s}\rangle )=k\leq \Cdim (A)$ where
$\d_i=\d_{i,n}+\d_{i,s}$ is the Jordan decomposition for $\d_i$.
\item If, an addition, $A$ is a central algebra then
$$ \GK (K\langle a_{1,
s}, \ldots , a_{t, s}\rangle )=k\leq \GK (A)\left(
1-\frac{1}{f_A+\max \{ f_A, 1\}} \right)$$ where $\d_{i,s}=\ad
(a_{i,s})$ for some $a_{i,s}\in A$, $f_A:= \max \{ \fdim (Q_m\t
A)\, | \, 0\leq m\leq  n\}$, and $\fdim $  is the (left) filter
dimension of the $Q_m$-algebra $Q_m \t A$.
\end{enumerate}
\end{corollary}

{\it Proof}. $1$. For each $i$, let $\d_i=\d_{i,n}+\d_{i,s}$ be
the Jordan decomposition for the locally finite derivation $\d_i$.
The derivations $\d_1, \ldots . \d_t$ commute, so $\D_s:=\{ \d_{1,
s},\ldots , \d_{t,s}\}$ is the set of commuting locally
semi-simple derivations of the algebra $A$ such that $\Ev (\D
)=\Ev (\D _s)$. So, without loss of generality one can assume that
all the derivations $\d_i$ are locally semi-simple.

By Corollary \ref{cGCgr3}, $E:=\Ev (\D )= \mathbb{Z}\l^1+\cdots +
\mathbb{Z}\l^k\subseteq K^t$ is a free abelian group of rank $k$
where $\l^1=(\l_i^1), \ldots , \l^k=( \l_i^k)$  are free
generators. Up to re-ordering  of the derivations $\d_1, \ldots ,
\d_t$ we may assume that the $k\times k$ matrix $\L =(\l^i_j)$,
$i,j=1,\ldots , k$, is nonsingular. Note that  $A=\bigoplus_{m\in
\mathbb{Z}^k}D(\D , m_1\l^1+\cdots +m_k\l^k)$ where $m=(m_1,\ldots
, m_k)$. For each $i=1, \ldots , k$, let us define a $K$-linear
map $\der_i:A\ra A$ that respects the $ \mathbb{Z}^k$-grading of
the algebra $A$ and acts in each space $D(\D , m_1\l^1+\cdots +
m_k\l^k)$ by multiplication on the scalar $\sum_{j=1}^km_j\l^j_i$.
By the very definition, all the maps  $\der_1, \ldots , \der_k$
commute and are locally semi-simple derivations of the algebra
$A$. Since all the derivations of the algebra $A$ are inner,
$\der_i =\ad (x_i)$ for some element $x_i\in A$. For each pair
$i\neq j$, $0=[\ad (x_i), \ad (x_j)]=\ad ([x_i, x_j])$, therefore
$\l_{ij}:=[x_i, x_j]\in Z(A)$, and  so $\l _{ij}=0$ since $\ad
(x_i)$ are locally semi-simple derivations. So, the elements $x_1,
\ldots , x_k$ commute. Let us show that they are algebraically
independent. Suppose that $f(x_1, \ldots , x_k)=0$  for a
polynomial $f(t_1, \ldots , t_k)\in K[t_1, \ldots , t_k]$. For
each nonzero element $a\in D(\D , \sum_{j=1}^km_j\l^j)$,
$0=af(x_1, \ldots , x_k)=f(x_1-\sum_{j=1}^km_j\l^j_1, \ldots ,
x_k-\sum_{j=1}^km_j\l^j_k)a$. So, $$f(x_1-\sum_{j=1}^km_j\l^j_1,
\ldots , x_k-\sum_{j=1}^km_j\l^j_k)=0, \;\;{\rm for \; all}\;\;
(m_1, \ldots , m_k)\in \mathbb{Z}^k.$$ This is possible iff $f=0$
since the $k\times k$ matrix $\L $ is non-singular and the  field
$K$ has characteristic zero. Then, $k\leq \Cdim (A)$.

Each derivation $\d_i$ is a locally semi-simple which acts on
$D(\D , m_1\l^1+\cdots +m_k\l^k)$ by multiplication on the scalar
$m_1\l^1_i+\cdots +m_k\l^k_i$. So, the Gelfand-Kirillov dimension
of the commutative subalgebra $K\langle \d_1, \ldots , \d_t\rangle
$ of ${\rm End}_K(A)$ is equal to the rank of the matrix
$(\l^j_i)$, that is $k$.

$2$. The second statement follows from statement 1 and its proof,
Theorem \ref{GKcsuba}, and the fact that the elements $a_{1,
s},\ldots , a_{t,s}$ commute. $\Box $


\section{Maximal Isotropic Subalgebras of Poisson Algebras}\label{maxissPA}

In this section, we apply Theorem \ref{GKcsuba} to obtain an upper
bound for the Gelfand-Kirillov dimension of (maximal) {\em
isotropic} subalgebras of certain Poisson algebras (Theorem
\ref{PoGKcsu}).

Let $(P, \{ \cdot , \cdot \} )$ be a {\em Poisson algebra} over
the field $K$. Recall that $P$ is an associative  commutative
$K$-algebra which is a Lie algebra with respect to the bracket $\{
\cdot , \cdot \}$ for which {\em Leibniz's rule} holds:
$$ \{ a, xy\} = \{ a, x\} y +x\{ a, y\}\;\; {\rm for \; all}\;\;
a,x,y\in P,$$
 which means that the {\em inner derivation } $\ad (a): P\ra P$,
 $x\mapsto \{ a, x\}$, of the Lie algebra $P$ is also a derivation
 of the associative algebra $P$. Therefore, to each Poisson
 algebra $P$ one can attach an associative subalgebra $A(P)$ of
 the ring of differential operators $\CD (P)$ with coefficients
 from the algebra $P$ which is generated by $P$ and $\ad (P):=\{
 \ad (a)\, | \, a\in P\}$. If $P$ is a finitely generated algebra
 then so is the algebra $A(P)$ with $\GK (A(P))\leq \GK (\CD
 (P))<\infty $.

 {\it Example}. Let $P_{2n}=K[x_1, \ldots , x_{2n}]$ be the {\em
 Poisson polynomial algebra} over a field $K$ of characteristic
 zero equipped with the {\em Poisson bracket}
 $$ \{ f,g\} =\sum_{i=1}^n( \frac{\der f}{\der x_i}\frac{\der g}{\der
 x_{n+i}}-\frac{\der f}{\der x_{n+i}}\frac{\der g}{\der
 x_i}).$$
 The algebra $A(P_{2n})$ is generated by the elements
 $$x_1, \ldots , x_{2n}, \; \ad (x_i)=\frac{\der }{\der
 x_{n+i}},\;
 \ad (x_{n+i})=-\frac{\der }{\der
 x_{i}}, \; i=1,\ldots , n.$$ So, the algebra $A(P_{2n})$ is canonically
 isomorphic to the Weyl algebra $A_{2n}$.

 Recall that the Weyl  algebra $A_n$ is the ring of differential
 operators $\CD ( \mathbb{A}^n)$ on the affine variety
 $\mathbb{A}^n$. As an abstract algebra the Weyl algebra $A_n$ is
 generated by $2n$ generators $x_1, \ldots, x_n,  \der_1, \ldots ,
 \der_n$ subject to the defining relations:
 $$ x_ix_j=x_jx_i, \; \der_i\der_j=\der_j\der_i, \;
 \der_ix_j-x_j\der_i=\d_{i,j}, \;\;  {\rm the \; Kronecker\;  delta},
 $$ for all $i,j=1, \ldots , n$.  The Weyl algebra $A_n$ is a
 central simple algebra of Gelfand-Kirillov dimension $2n$.

 {\it Definition}. We say that a Poisson algebra $P$ is a {\em
 strongly simple Poisson algebra} if
\begin{enumerate}
\item $P$ is a finitely generated (associative) algebra which is a
domain,
\item the algebra $A(P)$ is central simple, and
\item for each set of algebraically independent elements $a_1,
\ldots , a_m$ of the algebra $P$ such that $\{ a_i, a_j\}=0$ for
all $i,j=1, \ldots , m$ the (commuting) elements $a_1, \ldots ,
a_m, \ad (a_1), \ldots ,$ $ \ad (a_m)$ of the algebra $A(P)$ are
algebraically independent.
\end{enumerate}

\begin{theorem}\label{PoGKcsu}
Let $P$ be a strongly simple Poisson algebra, and $C$ be an
isotropic subalgebra of $P$, i.e. $\{ C,C\}=0$. Then
$$ \GK (C)\leq \frac{\GK (A(P))}{2}\left( 1-\frac{1}{f_{A(P)}+\max \{
f_{A(P)}, 1\}} \right)$$ where $f_{A(P)}:= \max \{
\fdim_{Q_m}(Q_m\t A(P))\, | \, 0\leq m\leq  \GK (A(P))\}$.
\end{theorem}

{\it Proof}. By the assumption the finitely generated algebra $P$
is a domain, hence the finitely generated algebra $A(P)$ is a
domain (as a subalgebra of the domain $\CD (Q(P))$, the ring of
differential operators with coefficients from the field of
fractions $Q(P)$ for the algebra $P$). It suffices to prove the
inequality for isotropic subalgebras of the Poisson algebra $P$
that are polynomial algebras. So, let $C$ be an isotropic
polynomial subalgebra of $P$ in $m$ variables, say  $a_1, \ldots ,
a_m$. By the assumption, the commuting elements $a_1,\ldots , a_m
, \ad (a_1), \ldots , \ad (a_m)$ of the algebra $A(P)$ are
algebraically independent. So,  the Gelfand-Kirillov dimension of
the subalgebra $C'$ of $A(P)$  generated by these elements is
equal to $2m$. By Theorem \ref{GKcsuba},
$$ 2 \GK (C)=2m=\GK (C')\leq \GK (A(P))\left( 1-\frac{1}{f_{A(P)}+\max \{
f_{A(P)}, 1\}} \right), $$ and this proves the inequality. $\Box $

\begin{corollary}\label{cPoGKcs}
\begin{enumerate}
\item The Poisson polynomial algebra $P_{2n}= K[x_1, \ldots ,
x_{2n}]$ (with the Poisson bracket)  over a field $K$ of
characteristic zero is a strongly simple Poisson algebra, the
algebra $A(P_{2n})$ is canonically isomorphic to the Weyl algebra
$A_{2n}$.
\item The Gelfand-Kirillov dimension of every  isotropic subalgebra of
the polynomial Poisson algebra $P_{2n}$ is $\leq n$.
\end{enumerate}
\end{corollary}

{\it Proof}. $1$. The third condition in the definition of
strongly simple Poisson algebra is the only statement we have to
prove. So, let $a_1, \ldots , a_m$ be  algebraically independent
elements of the algebra $P_{2n}$ such that $\{ a_i, a_j\}=0$ for
all $i,j=1, \ldots , m$. One can find polynomials, say
$a_{m+1},\ldots , a_{2n}$,  in $P_{2n}$ such that the elements
$a_1, \ldots , a_{2n}$ are algebraically independent, hence the
determinant $d$ of the Jacobian matrix $J:=(\frac{\der a_i}{\der
x_j })$ is a nonzero polynomial. Let $X=(\{ x_i, x_j\} )$ and
$Y=(\{ a_i, a_j\} )$ be, so-called, the {\em Poisson matrices}
associated with the elements $\{ x_i\}$ and $\{ a_i \} $. It
follows from $Y=J^TXJ$ that $\det (Y)=d^2 \det (X)\neq 0$ since
$\det (X)\neq 0$. The derivations
$$ \d_i:= d^{-1} \det
\begin{pmatrix}
  \{ a_1, a_1\}  & \ldots & \{ a_1 , a_{i-1} \}  & \{ a_1 , \cdot  \}  & \{ a_1 , a_{i+1} \}
   & \ldots  & \{ a_1 , a_{2n} \}  \\
\{ a_2, a_1\}  & \ldots & \{ a_2 , a_{i-1} \}  & \{ a_2 , \cdot
\} & \{ a_2 , a_{i+1} \}
   & \ldots  & \{ a_2 , a_{2n} \}  \\
   &  & & \ldots  &  & &  \\
\{ a_{2n}, a_1\}  & \ldots & \{ a_{2n} , a_{i-1} \}  & \{ a_{2n} ,
\cdot  \} & \{ a_1 , a_{i+1} \}
   & \ldots  & \{ a_{2n} , a_{2n} \}  \\
\end{pmatrix},
$$
$i=1, \ldots , 2n$, of the rational function  field $Q_{2n}=K(x_1,
\ldots , x_{2n})$ satisfy the following properties: $\d_i
(a_j)=\d_{i,j}$, the Kronecker delta. For each $i$ and $j$, the
kernel of  the derivation $\D_{ij}:=\d_i\d_j-\d_j\d_i\in
\Der_K(Q_{2n})$ contains $2n$ algebraically independent elements
 $a_1, \ldots , a_{2n}$. Hence $\D_{ij}=0$ since the field
 $Q_{2n}$ is algebraic over  its subfield $K(a_1, \ldots , a_{2n})$
 and ${\rm char} (K)=0$. So, the subalgebra, say $W$, of the ring
 of differential operators
 $\CD (Q_{2n})$ generated by the elements $a_1, \ldots , a_{2n},
\d_1, \ldots , \d_{2n}$ is isomorphic to the Weyl algebra
$A_{2n}$, and so $\GK (W)=\GK (A_{2n})=4n$.

Let $U$ be the $K$-subalgebra of $\CD (Q_{2n})$ generated by the
elements $x_1, \ldots , x_{2n}, \d_1, \ldots , \d_{2n}$, and
$d^{-1}$. Let $P'$ be the localization of the polynomial algebra
$P_{2n}$ at the powers of the element $d$. Then $\d_1, \ldots ,
\d_{2n}\in \sum_{i=1}^{2n}P'\ad (a_i)$ and $\ad (a_1), \ldots ,
\ad (a_{2n})\in  \sum_{i=1}^{2n}P'\d_i$, hence the algebra $U$ is
generated (over $K$) by $P'$ and  $\ad (a_1) , \ldots , \ad
(a_{2n})$. The algebra $U$ can be viewed as  a subalgebra of the
ring of differential operators $\CD (P')$. Now, the inclusions,
$W\subseteq U\subseteq \CD (P')$ imply $4n=\GK (W)\leq \GK (U)\leq
\GK (\CD (P'))=2\GK (P')=4n$, therefore $\GK (U)=4n$. The algebra
$U$ is a factor algebra of an iterated Ore extension $V=P'[t_1;\ad
(a_1)]\cdots [t_{2n}; \ad (a_{2n})]$. Since $P'$ is a domain, so
is the algebra $V$. The algebra $P'$ is a finitely generated
algebra of Gelfand-Kirillov dimension $2n$, hence $\GK (V)=\GK
(P')+2n =4n$ (by \cite{MR}, 8.2.11). Since $\GK (V)=\GK (U)$ and
any proper factor algebra of $V$ has Gelfand-Kirillov dimension
strictly less than $\GK (V)$ (by \cite{MR}, 8.3.5, since $V$ is a
domain), the algebras $V$ and $U$ must be isomorphic. Therefore,
 the (commuting) elements $a_1, \ldots
, a_m, \ad (a_1), \ldots ,$ $ \ad (a_m)$ of the algebra $U$ (and
$A(P)$) must be  algebraically independent.

$2$. Let $C$ be an isotropic subalgebra of the Poisson algebra
$P_{2n}$. Note that $f_{A(P_{2n})}=f_{A_{2n}}=1$  and $\GK
(A_{2n})=4n$. By Theorem \ref{PoGKcsu},
$$ \GK (C) \leq \frac{4n}{2}(1-\frac{1}{1+1})=n. \; \Box $$

{\it Remark}. This result means that for the Poisson polynomial
algebra $P_{2n}$  the right hand side in the inequality of Theorem
\ref{PoGKcsu} is the {\em exact} upper bound for the
Gelfand-Kirillov dimension of isotropic subalgebras in $P_{2n}$
since the polynomial subalgebra $K[x_1, \ldots, x_n]$ of $P_{2n}$
is isotropic.


\section{Holonomic Modules}\label{sholmod}

{\it Definition}. Let $A$ be a finitely generated $K$-algebra, and
$h_A$ be its holonomic number. A nonzero finitely generated
$A$-module $M$ is called a {\em holonomic} $A$-module if $\GK
(M)=h_A$. We denote  by $\hol (A)$ the set of all the holonomic
$A$-modules.

Since the holonomic number is an infimum it is not clear at the
outset that there will be modules which achieve this dimension.
Clearly, $\hol(A)\neq \emptyset $ if the Gelfand-Kirillov
dimension of every finitely generated $A$-module is a natural
number.

A nonzero submodule or a factor module of a holonomic is a
holonomic module (since the Gelfand-Kirillov dimension of a
submodule or a factor module does not exceed the Gelfand-Kirillov
of the module). If, in addition, the finitely generated algebra
$A$ is left Noetherian and finitely partitive then each holonomic
$A$-module $M$ has finite length and each simple sub-factor of $M$
is a holonomic module.

\begin{lemma}\label{AMBAM}
Let $A$ and $B$ be finitely generated  $K$-algebras, and ${}_AM_B$
be a bimodule such that ${}_AM$ is finitely generated. Then $\GK
({}_AM_B)\leq \GK ({}_AM)$.
\end{lemma}

{\it Proof}. 
Let $M_0$ be a finite dimensional generating subspace for the
$A$-module $M$, and let $\{ A_i\}$ and $\{ B_i\}$ be standard
(finite dimensional) filtrations for the algebras $A$ and $B$
respectively. Then $M_0B_1\subseteq A_nM_0$ for some $n\geq 0$.
Now, $\{ M_i:=\sum_{j=0}^iA_jM_0B_{i-j}\}$ is the standard finite
dimensional filtration for the bimodule ${}_AM_B$. Obviously,
$$M_i=\sum_{j=0}^iA_jM_0B_1^{i-j}\subseteq
\sum_{j=0}^iA_jA_{n(i-j)}M_0\subseteq A_{i(n+1)}M_0\;\;  {\rm
for\;   all}\;\;  i\geq 0.$$
 Hence, $\GK ({}_AM_B)\leq \GK
({}_AM)$.

$\Box $

\begin{theorem}\label{simhol}
Let a finitely generated $K$-algebra $A$ be a domain with $0<\GK
(A)<\infty$. Suppose that $C$ is a commutative finitely generated
subalgebra of $A$ with field of fractions $Q$ such that $\GK
(A)-\GK (C)=h_{A\t Q}$, the holonomic number for the $Q$-algebra
$A\t Q$. Then $A\t_CQ$ is a simple holonomic module over the
$Q$-algebra $A\t Q$ (i.e. $\GK_Q({}_{A\t Q}A\t_CQ)=h_{A\t Q}$).
\end{theorem}

{\it Proof}. Since $\GK (C)\leq \GK (A)$, the holonomic number
$h_{A\t Q}=\GK (A)-\GK (C)<\infty$.

The $A\t Q$-module $A\t_CQ$ is a nonzero module. By Proposition
\ref{drGK},
$$ \GK (A)=\GK ({}_AA_A)\geq \GK ({}_AA_C)=\GK ({}_{A\t C}A)\geq
\GK_Q({}_{A\t Q}(A\t_CQ))+\GK (C),$$ hence
$$\GK_Q({}_{A\t Q}(A\t_CQ))\leq \GK (A)-\GK (C)=h_{A\t Q}. $$
This means that $A\t_CQ$ is a holonomic module  of the $Q$-algebra
$A\t Q$.

The quotient field $Q$ for the algebra $C$ is the localization
$C\S1 $ of the domain $C$ at its multiplicatively closed subset
$S:=C\backslash \{ 0\}$. So, $A\t_CQ\simeq A\S1 $ is the right
localization of the right $C$-module $A$ at $S$, and the left
localization of the left $A\t C$-module $A$ (i.e. ${}_{A\t
C}A={}_AA_C$) at $S$ considered as the subset $\{ 1\t c\, | \,
c\in S\}$ of $A\t C$. The algebra $A\t Q$ is a localization of the
algebra $A\t C$ at $S$.  Since $A$ is a domain and $S\subseteq A$,
the natural map $A\ra A\t_CQ\simeq A\S1$ is an $A\t C$-module
monomorphism. So, we identify $A$ in $A\S1$. Suppose that $A\t_CQ$
is not a simple $A\t Q$-module. Then one can find a nonzero proper
$A\t Q$-submodule, say $M$, of $A\t_CQ$ (i.e. $0\neq M\neq
A\t_CQ$). We seek a contradiction. Then $N:=A\cap M$ is a nonzero
$A\t C$-module since $M=N\S1$.

Localizing the short exact sequence of $A\t C$-modules: $0\ra N\ra
A\ra A/N\ra 0$ at $S$ we get a short exact sequence of $A\t
Q$-modules:
$$0\ra M\ra A\S1 \ra L:=(A/N)\S1 \ra 0,$$ with $L\neq 0$
since $M\neq A\S1 $. Fix an arbitrary  nonzero element, say $a$ of
$N$. The algebra $A$ is a domain, so the $A$-submodule $Aa$ of $N$
is isomorphic to ${}_AA$. By \cite{MR}, 8.3.5,
$$\GK
({}_A(A/Aa))\leq \GK ({}_AA)-1<\GK (A).$$
 The $A$-module $A/N$ is
an epimorphic image of the $A$-module $A/Aa$, hence
\begin{eqnarray*}
\GK (A)&>& \GK ({}_A(A/Aa))\geq \GK ({}_A(A/N))\\
&\geq & \GK ({}_{A\t C}(A/N))\;\;\; ({\rm by}\;\; {\rm Lemma}\;
\ref{AMBAM})\\
&\geq & \GK_Q ({}_{A\t Q}L)+\GK (C)\;\;\; ({\rm by\; Proposition
}\; \ref{drGK}).
\end{eqnarray*}
Now,
$$ h_{A\t Q}\leq \GK_Q({}_{A\t Q}L)<\GK (A)-\GK (C)=h_{A\t Q},
\;\; {\rm a \; contradiction}.$$ So, the $A\t Q$-module $A\t_CQ$
must be simple. $\Box $

\begin{corollary}\label{csimhol}
 Let $K$ be an algebraically closed  field of characteristic zero,
$X$ be a smooth irreducible affine  algebraic variety of dimension
$n:=\dim (X)>0$, and $C$ be a commutative subalgebra of the ring
of differential operators $\CD (X)$ on $X$ with $\GK (C)=n$, $Q$
be the field of fractions for $C$. Then $\CD (X)\t_CQ$ is a simple
holonomic module over the $Q$-algebra $\CD (X)\t Q$ (i.e.
$\GK_Q({}_{\CD (X)\t Q}\CD (X)\t_C Q=n$).
\end{corollary}

{\it Proof}. Since $\GK (\CDX (X))=2n$ and $h_{\CD (X)\t Q}=n$,
the result follows from Theorem \ref{simhol}. $\Box $

Department of Pure Mathematics

University of  Sheffield

Hicks Building

Sheffield S3~7RH

UK

email: v.bavula@sheffield.ac.uk




\begin{thebibliography}{99}

\bibitem{AmitsurPAMS58} S. A. Amitsur,  Commutative linear differential operators.
{\it  Pacific J. Math.} {\bf  8} (1958) 1--10.

\bibitem{Amitsur-Small78} S. A. Amitsur and L. W. Small,  Polynomials over division rings.
{\it Israel J. Math.} {\bf  31} (1978), no. 3-4, 353--358.

\bibitem{Bavcafd} V. Bavula, Filter dimension of algebras and modules, a simplicity
criterion for generalized Weyl algebras. {\it Comm. Algebra} {\bf
24} (1996), no. 6, 1971--1992.

\bibitem{Bavjafd} V. Bavula, Krull, Gelfand-Kirillov, and filter dimensions of simple affine
   algebras. {\it J. Algebra} {\bf  206} (1998), no. 1, 33--39.

\bibitem{bie98} V. Bavula, Krull, Gelfand-Kirillov, filter, faithful and Schur dimensions.
   {\it Infinite length modules} (Bielefeld, 1998), 149--166, Trends Math., Birkhäuser,
   Basel, 2000.

\bibitem{Ber72} I.N. Bernstein,  Analytic continuation of generalized
 functions with respect to a parameter. Funkcional. Anal. i Prilozen.
  {\bf 6} (1972), no. 4, 26--40.


\bibitem{BrEtGaQCIS}  A. Braverman, P. Etingof, and D. Gaitsgory,  Quantum integrable systems and
   differential Galois theory. Transform. Groups, {\bf  2} (1997), no. 1, 31--56.


\bibitem{Dix}
J. Dixmier,  Sur les alg\`{e}bres de Weyl, {\it Bull. Soc. Math.
France} {\bf 96} (1968), 209--242.

\bibitem{DixUEA}
J. Dixmier,  {\em Alg\'{e}bres enveloppantes}, Gauthier-Villars,
Paris, 1974.


\bibitem{GK66} I. M. Gelfand and A. A. Kirillov,
Sur les corps li\'es aux alg\`ebres enveloppantes des alg\`ebres
de Lie, {\it Publ. Math. IHES}, {\bf  31} (1966), 5--19.

\bibitem{Joseph72HP} A. Joseph, Gelfand-Kirillov dimension for
algebras associated with the Weyl algebra, Ann. Inst. H.
Poincar\'e, {\bf 17} (1972), 325--336.

\bibitem{JosLN74} A. Joseph, Sur les alg\`ebres de Weyl, Lecture
 notes, 1974 (unpublished).

\bibitem{JosgenQ} A. Joseph,  A generalization of Quillen's lemma and its application to
the Weyl   algebras. {\it Israel J. Math.},  {\bf 28} (1977), no.
3, 177--192.

\bibitem{KL} G. Krause and T.  Lenagan, {\em Growth of algebras and Gelfand-Kirillov
 dimension}. Revised edition. Graduate Studies in Mathematics, 22.
American Mathematical Society, Providence, RI, 2000.

\bibitem{Mak-LimFree} L. Makar-Limanov,  On subalgebras of the first Weyl skewfield.
{\it Comm. Algebra},   {\bf 19} (1991), no. 7, 1971--1982.

\bibitem{Mak-Limcom} L. Makar-Limanov,  Commutativity of certain centralizers in the rings
   $R\sb{n,\,k}$. (Russian) {\it Funkcional. Anal. i Prilozen}. {\bf 4} (1970), no. 4, 78.

\bibitem{MR} J. C. McConnell and J. C. Robson, Noncommutative Noetherian rings,
Wiley, Chichester,  1987.

\bibitem{Pierceb}  R. Pierce, Associative algebras. Springer-Verlag, New York-Berlin, 1982.

 \bibitem{Ren-Gab} R. Rentschler and P. Gabriel, Sur la dimension des anneaux et ensembles
   ordonn\'es. (French) {\it C. R. Acad. Sci. Paris Sér. A-B} {\bf 265} (1967), A712--A715.

\bibitem{Resco79} R. Resco, Transcendental division algebras and simple
Noetherian rings.  {\it Israel J. Math.},  {\bf 32} (1979), no.
2-3, 236--256.


\end{thebibliography}
\end{document}